\documentclass{article}

\newtheorem{thm}{Theorem}[section]
\newtheorem{cor}[thm]{Corollary}

\newtheorem{lem}[thm]{Lemma}

\newtheorem{ex}[thm]{Example}

\newcommand{\be}{\begin{equation}}
\newcommand{\ee}{\end{equation}}
\newcommand{\ben}{\begin{enumerate}}
\newcommand{\een}{\end{enumerate}}
\newcommand{\pa}{{\partial}}
\newcommand{\R}{{\rm R}}
\newcommand{\e}{{\epsilon}}
\newcommand{\g}{{\bf g}}
\newcommand{\pxi}{{\pa \over \pa x^i}}
\newcommand{\pxj}{{\pa \over \pa x^j}}

\newcommand{\qed}{\hspace*{\fill}Q.E.D.}  
\title{\Large Finsler Manifolds with Nonpositive Flag Curvature\\
 and Constant S-curvature}
\author{Zhongmin Shen}

\date{}

\begin{document}

\maketitle

\begin{abstract} 
The flag curvature is a natural extension of the sectional curvature in Riemannian geometry, and the S-curvature  is a non-Riemannian quantity which  vanishes for Riemannian metrics. There are (incomplete) non-Riemannian  Finsler metrics on an open subset in $\R^n$ with negative flag curvature and constant S-curvature. In this paper, we are going to show a global rigidity theorem that every Finsler metric with negative flag curvature and constant S-curvature must be Riemannian if the manifold is compact. We also study the nonpositive flag curvature case.
\end{abstract}

\bigskip

\section{Introduction}

One of  important problems in Finsler geometry is to understand the geometric  meanings of various quantities and their impacts on the global geometric structures.  
 Imaging a Finsler manifold as an Easter egg and a Riemannian manifold as a white egg, Finsler manifolds are not only curved, but also very ``colorful''. The flag curvature ${\bf K}$ tells us how curved is the Finsler manifold at a point. 
There are several non-Riemannian quantities which describe the ``color'' and its rate of change over the manifold,   such as the mean Cartan torsion ${\bf I}$, the mean Landsberg curvature ${\bf J}$
and the S-curvature  ${\bf S}$ (see \cite{Sh2} or Section \ref{sectionP} below). These quantities interact with the flag curvature  in a delicate way. 
 The  mean Landsberg curvature  and the S-curvature reveal  different non-Riemannian properties. 
For examples, there is a family of Finsler metrics on $S^3$ with ${\bf K}=1$ and ${\bf S}=0$ \cite{BaSh}. However, every local Finsler metric with ${\bf K}=1$ and ${\bf J}=0$ must be Riemannian (Theorem 9.1.1 in \cite{Sh2}).  

An $n$-dimensional  Finsler metric is said to have {\it constant S-curvature} if ${\bf S}= (n+1) c F$ for some constant $c$. It is known that every Randers metric of constant flag curvature has constant S-curvature \cite{BaRo}, \cite{BaRoSh}. 
This is one of our motivations  to consider Finsler metrics of constant S-curvature.  
In this paper, we are going to prove the following global metric rigidity theorem.

\begin{thm}\label{thmnonpositive} 
Let $(M, F)$ be an $n$-dimensional compact boundaryless Finsler manifold  with constant S-curvature, i.e., ${\bf S}= (n+1) c F$ for some constant $c$. 
\ben
\item[(a)] If $F$ has negative flag curvature, ${\bf K} <0$, then it must be Riemannian;
\item[(b)] If $F$ has nonpositive flag curvature, ${\bf K}\leq 0$, then the mean Landsberg curvature vanishes, ${\bf J}=0$,  and  the flag curvature ${\bf K}(P,y)=0$  for the flags $P= {\rm span}\{y, {\bf I}_y \}\subset T_xM$ whenever ${\bf I}_y\not=0$. 
\een
\end{thm}

The compactness in Theorem \ref{thmnonpositive} (a) can not be dropped. 
Consider the following family of Finsler metrics on the unit ball ${\rm B}^n \subset \R^n$,
\be
F = { \sqrt{ |y|^2 - (|x|^2|y|^2 -\langle x, y \rangle^2}+\langle x, y \rangle \over 1-|x|^2 } + { \langle a, y \rangle \over 1+\langle a, x \rangle },  \label{FunkFa}
\ee
where $y\in T_x{\rm B}^n \cong \R^n$ and $a\in \R^n$ is an arbitrary constant vector with $|a| < 1$. 
 It is proved that 
$F$ has constant flag curvature ${\bf K}=-{1\over 4}$ and constant S-curvature ${\bf S} ={1\over 2} (n+1) F$  (see \cite{Sh2}\cite{Sh3}). Clearly, $F$  is not Riemannian.

The compactness in Theorem \ref{thmnonpositive} (b) can not be  dropped. 
Let $n\geq 2$ and
\[ {\cal U}:=\Big \{ p=(s,t, \bar{p}) \in \R^2 \times \R^{n-2} 
\ \Big | \ s^2+t^2 < 1 \Big \} .\]
Define
\[
 F:= {\sqrt{ \Big ( -t u + s v \Big )^2 
+ |y|^2 \Big ( 1- s^2 -t^2 \Big ) } - \Big ( - tu + s v \Big ) 
\over 1- s^2 - t^2 },
\]
where $y = (u, v,\bar{y} ) \in T_p {\cal U} \cong \R^n$ and $p= (s,t, \bar{p})\in {\cal U}$.
$F$ is an incomplete Finsler metric on $\Omega$ with ${\bf K}=0$ and ${\bf S}=0$, but ${\bf J}\not=0$ \cite{Sh5}.

The compactness condition in Theorem \ref{thmnonpositive} can be replaced by a completeness condition together with certain growth condition on the mean Cartan torsion. See Theorems \ref{thmnonpositive*} and \ref{thm4.2} below. 

\bigskip

\begin{cor}
Let $(M, F)$ be a compact boundaryless Berwald manifold with nonpositive flag curvature. Then the following hold, 
\ben
\item[(a)] If $F$ has negative flag curvature, ${\bf K} <0$, then it must be Riemannian;
\item[(b)] If $F$ has nonpositive flag curvature, ${\bf K}\leq 0$, then   ${\bf K}(P,y)=0$ for the flag $P={\rm span}\{ y, {\bf I}_y \}$ whenever ${\bf I}_y\not=0$.
\een
\end{cor}

In dimension two, we have the following
\begin{cor}
Let $(M, F)$ be a compact boundaryless Finsler surface. Suppose that   ${\bf K} \leq 0$ and  ${\bf S}=3 c F$ for some constant $c$,  then $F$ is either locally Minkowskian or Riemannian.
\end{cor}
The proof is simple. First,  by  Theorem \ref{thmnonpositive*} below, we know that  ${\bf J}=0$, then  the theorem follows from Theorem 7.3.2 in \cite{BaChSh2}.

\bigskip

In dimension $n\geq 3$, 
we  have some non-trivial examples satisfying the conditions and conclusions in Theorem \ref{thmnonpositive} (b). Let $(N, h)$ be an arbitrary closed hyperbolic Riemannian manifold. For any $\e \geq 0$, let 
\[ F_{\e}:= \sqrt{ h^2(\bar{x}, \bar{y}) + w^2 +\e \sqrt{ h^4 (\bar{x}, \bar{y}) + w^4 } },\]
where $ x =(\bar{x}, s)\in M$ and $ y = \bar{y}\oplus w{\pa \over \pa s}\in T_xM$. 
This family of Finsler metrics is constructed by Z.I. Szab\'{o} in his classification of Berwald metrics \cite{Sz}. It is known that  each $F_{\e}$ is a Berwald metric. Thus ${\bf J}=0$ and ${\bf S}=0$ \cite{Sh2}. Further it can be shown  that $F_{\e}$ satisfies
 that ${\bf K} \leq 0$ and ${\bf K}(P, y)=0$ for $P={\rm span}\{y,  {\bf I}_y \}$. The proof will be given in Section \ref{sectionExample} below. 
A natural problem arises: Is the Finsler metric in Theorem \ref{thmnonpositive} (b)  a Berwald metric? This problem remains open.

Finally, we should point out that there are already
several global rigidity results on the metric structure of Finsler manifolds with ${\bf K} \leq 0$. For example, H. Akbar-Zadeh proves that every closed Finsler manifold with ${\bf K}= -1$ must be Riemannian and every closed Finsler manifold with ${\bf K}=0$ must be locally Minkowskian  \cite{AZ}. Mo-Shen prove that every closed Finsler manifold of scalar curvature with ${\bf K}<0$ must be of Randers type in dimension $\geq 3$ \cite{MoSh}. Here a Finsler metric
$F$ is said to be {\it of scalar curvature} if the flag curvature ${\bf K}={\bf K}(x,y)$ is independent of $P$ for any given direction $y\in T_xM$. Riemannian metrics of scalar curvature must have isotropic sectional curvature ${\bf K}={\bf K}(x)$, hence they have  constant sectional curvature in dimension $n\geq 3$ by the Schur Lemma. But there are lots of Finsler metrics of scalar curvature which have not been completely classified yet.

\section{Preliminaries}\label{sectionP}

In this section, we are going to give a brief description on the flag curvature and the above mentioned  non-Riemannian quantities.

Let $M$ be an $n$-dimensional manifold and let
$\pi: TM_o:=TM\setminus\{0\}\to M$
 denote the slit tangent bundle. The pull-back tangent bundle is defined by
 $\pi^*TM := \{ (x,y, v) \ | \ 0\not=y, v\in T_xM\}$ and the pull-back cotangent bundle is defined by $\pi^*T^*M:= \{ \pi^*\theta \ | \ \theta\in T^*M\}$. 

By definition, a Finsler metric $F$ on a manifold $M$ is a nonnegative function on $TM$ which is positively $y$-homogeneous  of degree one with positive definite fundamental tensor
 ${\bf g} := g_{ij}dx^i \otimes dx^j $ on $\pi^*TM$, where $g_{ij} := {1\over 2} [F^2]_{y^iy^j}(x,y)$. A  special class of Finsler metrics  are Randers metrics in the form $F=\alpha+\beta$ where $\alpha=\sqrt{a_{ij}(x)y^iy^j}$ is a Riemannian metric and $\beta = b_i(x)y^i$ is a $1$-form  with $ \|\beta \|_{x} :=\sqrt{ a^{ij}(x)b_i(x)b_j(x)} <1$ for any $x\in M$.

For a Finsler metric $F$, the volume $dV = \sigma_F (x) dx^1 \cdots dx^n$  is defined by 
\be
 \sigma_F(x) := {{\rm Vol} ({\rm B}^n(1))\over 
 {\rm Vol} \Big \{ (y^i) \in \R^n \Big | \ F\Big (x, \; y^i \pxi|_x \Big ) < 1 \Big \} }.\label{sigma_F}
\ee
When $F=\sqrt{g_{ij}(x)y^iy^j}$ is Riemannian, then $\sigma_F(x)=\sqrt{\det(g_{ij}(x))}$. 
In general, the following quantity is not equal to zero,
\[ \tau(x,y) :=\ln  \Big [ { \sqrt{\det ( g_{ij}(x,y) )}\over \sigma_F(x)} \Big ] .\]
$\tau=\tau(x,y)$ is a scalar function on $TM_o$, which  is called  the  {\it distortion} \cite{Sh2}. The  distortion is our primary non-Riemannian quantity. 
Let 
\be
 I_i:= {\pa \tau \over \pa y^i}(x,y) = {1\over 2} g^{jk}(x,y) {\pa g_{jk}\over \pa y^i}(x,y) .\label{I_i}
\ee
We have
\be
I_i y^i = 0. \label{I=0}
\ee 
The tensor ${\bf I}:= I_i dx^i$ on $TM_o$  is called the {\it mean Cartan tensor}.  
According to Deicke's theorem \cite{De}, $F$ is Riemannian if and only if ${\bf I}=0$. 
Define 
the norm of ${\bf I}$ at a point $x\in M$ by
\[ \| {\bf I} \|_x :=\sup_{0\not=y \in T_xM} \sqrt{ I_i (x, y) g^{ij} (x, y) I_j (x, y) }.\] 
For a point $p\in M$, let 
\[ {\cal I}_p(r):= \sup_{\min (d(p, x), d(x,p)) < r } \| {\bf I} \|_x.\]
The mean Cartan tensor ${\bf I}$ is said to grow sub-linearly if for any point $p\in M$,
\[  {\cal I}_p(r) = o(r), \ \ \ \ \ \ \ (r \to +\infty).\]
${\bf I}$ is said to grow {\it sub-exponentially} at rate of $k=1$ if for any point $p\in M$,

\[ {\cal I}_p(r) = o(e^r), \ \ \ \ \ \ \ (r\to +\infty).\]

 It is known that for a Randers metric $F=\alpha+\beta$, ${\bf I}$ is bounded, i.e.,
\[
\| {\bf I}\|_x \leq {n+1\over \sqrt{2}} \sqrt{ 1 - \sqrt{1-\|\beta \|_x^2 }} < {n+1\over \sqrt{2}}, \ \ \ \ \ x\in M.\]
The bound in dimension two is suggested by B. Lackey. 
See Proposition 7.1.2 in \cite{Sh2} for a proof. 

The geodesics in a Finsler manifold are characterized by a system of second order ordinary differential equations
\[ \ddot{\sigma}^i + 2G^i (\sigma, \dot{\sigma}) =0,\]
where $G^i = G^i(x,y)$ are positively $y$-homogeneous functions of degree two. When $F$ is Riemannian, $G^i = {1\over 2}\Gamma^i_{jk}(x)y^jy^k$ are quadratic in $y\in T_xM$. A Finsler metric with such a property  called a {\it Berwald metric}.   There are  many non-Riemannian Berwald manifolds (see Section \ref{sectionExample} below).

For a non-zero vector $y\in T_xM$, set
\[ {\bf S}(x,y) := {d\over dt} \Big [ \tau \Big ( \sigma(t), \dot{\sigma}(t) \Big ) \Big ]|_{t=0},\]
where $\sigma=\sigma(t)$ is the geodesic with $\sigma(0)=x$ and $\dot{\sigma}(0)=y$. 
 ${\bf S}={\bf S}(x,y)$ is a scalar function  on $TM_o$ which is called the {\it S-curvature} \cite{Sh1}\cite{Sh2}. 
Let $dV = \sigma_F(x) dx^1 \cdots dx^n$ be the  volume form on $M$. 
The S-curvature can be expressed by 
\be {\bf S} = {\pa G^m \over \pa y^m}(x,y) - y^m {\pa \over \pa x^m} \Big [ \ln \sigma_F(x) \Big ]. \label{Scurvature}
\ee
 It is proved that ${\bf S}=0$ for Berwald metrics \cite{Sh1}\cite{Sh2}. 
An $n$-dimensional  Finsler metric $F$ is said to {\it have constant S-curvature} if there is a constant $c$ such that $ {\bf S}= (n+1) c F$. It is known that all Randers metrics of constant flag curvature must  have constant S-curvature \cite{BaRo} (see \cite{BaRoSh} for the classification of such metrics).

\bigskip

There is a distinguished linear connection  $\nabla$  on $\pi^*TM$ which is called the {\it Chern connection} \cite{Ch}. Let $\{ {\bf e}_i\}$ be a local frame for $\pi^*TM$ and $\{\omega^i\}$ the dual local frame for $\pi^*T^*M$. $\nabla$ can be expressed by 
\[ \nabla V  = \Big \{ dV^i + V^j \omega_j^{\ i}\Big  \} \otimes {\bf e}_i,\]
where $V= V^i {\bf e}_i\in C^{\infty}(\pi^*TM)$.  The Chern connection can be viewed as a generalization of the Levi-Civita connection in Riemannian geometry. 
Let 
\[ \omega^{n+i} := dy^i + y^j \omega_j^{\ i},\]
where $y^i$ are local functions on $TM_o$ defined by the canonical section ${\bf Y}= y^i {\bf e}_i$ of $\pi^*TM$.
We obtain a local coframe $\{ \omega^i, \omega^{n+i}\}$ for $T^*(TM_o)$.

Let
\[
 \Omega^i := d \omega^{n+i} - \omega^{n+j} \wedge \omega_j^{\ i} .\]
 $\Omega^i$ can be expressed as follows, 
\[ \Omega^i = {1\over 2} R^i_{\ kl} \omega^k \wedge \omega^l - L^i_{\ kl}\omega^k \wedge \omega^{n+l},\]
where $R^i_{\ kl}+ R^i_{\ lk}=0$ and $L^i_{\ kl}=L^i_{\ lk}$. 
The anti-symmetric tensor ${\bf R} = R^i_{\ kl}{\bf e}_i \otimes \omega^k \otimes \omega^l$ is called the {\it Riemann tensor} and the symmetric tensor
${\bf L} = L^i_{\ kl} {\bf e}_i \otimes \omega^k \otimes \omega^l$ is called the Landsberg tensor.

Let 
\[ R^i_{\ k} := R^i_{\ kl} y^l, \ \ \ \ \ \ R_{jk}:=g_{ij}R^i_{\ k}. \] 
We have
\be
 R^i_{\ k} y^k=0, \ \ \ \ \ R_{jk}= R_{kj}.
\label{R=0}
\ee
See \cite{Sh2} for details. 
The tensor ${\bf R}:= R^i_{\ k}{\bf e}_i \otimes \omega^k$  is still called the {\it Riemann tensor}.
The notion of Riemann (curvature) tensor for general Finsler metrics is introduced by L. Berwald using the Berwald connection \cite{Be1}\cite{Be2}. 
Let 
\[
J_k := L^m_{\ km}.\]
 The tensor
${\bf J}:= J_i \omega^i$   is called the {\it mean Landsberg tensor}. For a Berwald metric, ${\bf J}=0$ \cite{Sh2}.

For a scalar function on $TM_o$,  say  $\tau$, we define its covariant derivatives by  
\[
 d\tau = \tau_{|k} \omega^k + \tau_{\cdot k} \omega^{n+k}.\]
From (\ref{I_i}), we have 
\[
 \tau_{\cdot i} = {\pa \tau \over \pa y^i} = I_i .
\]
We have
\[
 {\bf S} := \tau_{|m} y^m.\]

For a tensor, say,  ${\bf I}= I_{i} \omega^{i}$,  the covariant derivatives are defined in a canonical way  by
\[
dI_{i} - I_k  \omega_{i}^{\  k}
  = I_{i|k} \omega^k + I_{i \cdot k} \omega^{n+k}.
\]
We have
\be
J_i  =
 I_{i|m} y^m  \label{Pg3***}
 \ee
Hence
\be
J_i y^i =0. \label{J=0}
\ee
See \cite{Sh2} for details. 

\bigskip

Now we interpret the above geometric quantities in a different way.

Let $F$ be a Finsler metric on an $n$-dimensional manifold $M$. 
For a non-zero tangent vector $ y=y^i\pxi|_x\in T_xM$, define 
\[ {\bf g}_y(u, v):= g_{ij}(x,y) u^i v^j, \ \ \ \ \ \ u=u^i\pxi|_x, v =v^j\pxj|_x\in T_xM,\]
where $g_{ij}(x,y) = {1\over 2} [F^2]_{y^iy^j}(x,y)$.
 Each ${\bf g}_y$ is an inner product 
on the tangent space $T_xM$.

The  Riemann tensor can be viewed as a family of endomorphisms on tangent spaces. 
\[ {\bf R}_y (u): = R^i_{\ k}(x,y) u^k \pxi|_x,
\]
where $u=u^i\pxi|_x \in T_xM$. The coefficients $R^i_{\ k}=R^i_{\ k}(x,y)$ are given by
\be
R^i_{\ k} = 2 {\pa G^i\over \pa x^k}-y^j{\pa^2 G^i\over \pa x^j\pa y^k}
+2G^j {\pa^2 G^i \over \pa y^j \pa y^k} - {\pa G^i \over \pa y^j}{\pa G^j \over \pa y^k}.  \label{Riemann}
\ee
 It follows from (\ref{R=0})  that
\be  {\bf R}_y (y)=0, \ \ \ \ {\bf g}_y ({\bf R}_y (u), v)= {\bf g}_y (u, {\bf R}_y (v) ),
\label{IJK}
\ee
where $u, v\in T_xM$.
The family ${\bf R}:=\{ {\bf R}_y | y\in T_xM\setminus\{0\}\}$ is called  the {\it Riemann curvature}.

Using the Chern connection $\nabla$ on $\pi^*TM$, one can define the covariant derivative of  
a vector field $X= X^i (t) \pxi |_{c(t)} $ along a curve $c$ by 
\[ D_{\dot{c}} X(t):= \Big \{ {d X^i \over dt}(t) +  X^j(t) \Gamma^i_{jk}( c(t), \dot{c}(t)) \dot{c}^k(t) \Big \} \pxi |_{c(t)}.\]
If $H=H(s,t)$ is a family of geodesics, i.e., for each $s$, $\sigma_s(t):= H(s,t)$ is a geodesic, the variation field $V_s(t):= {\pa H\over \pa s}(s, t)$ satisfies
the following Jacobi field along $\sigma_s$,
\[ D_{\dot{\sigma}_s} D_{\dot{\sigma}_s} V_s(t) + {\bf R}_{\dot{\sigma}_s(t)} (V(t)) =0.\]

For a tangent plane $P \subset T_xM$ and a vector $0\not=y\in P$, let
\[ {\bf K}(P, y):= { {\bf g}_y ({\bf R}_y(u), u) )\over 
{\bf g}_y (y, y) {\bf g}_y(u, u) - [ {\bf g}_y (y, u) ]^2 } ,\]
where $P= {\rm span}\{ y, u\}$. By (\ref{IJK}), one can see that ${\bf K}(P, y)$ is well-defined, namely, independent of the choice of a particular $u\in T_xM$. 

The mean Cartan tensor and the mean Landsberg tensor can be viewed as families of vectors on the manifold, i.e., 
\[ {\bf I}_y = I^i (x, y) \pxi|_x, \ \ \ \ \ \ {\bf J}_y = J^i (x, y) \pxi|_x,\]
where  $I^i:= g^{il} I_l$ and $J^i:= g^{il} J_l$.
It follows from (\ref{I=0}) and (\ref{J=0}) that
\[ {\bf g}_y ({\bf I}_y, y) = 0 = {\bf g}_y ({\bf J}_y, y).\] 
Thus ${\bf I}_y$ and ${\bf J}_y$ are perpendicular to $y$ with respect to ${\bf g}_y$. 
We call ${\bf I}:=\{ {\bf I}_y \ | \  y\in TM\setminus\{0\}\}$ 
and ${\bf J}:=\{ {\bf J}_y\  | \ y\in TM\setminus\{0\}\}$ the {\it mean Cartan torsion}
and the {\it Landsberg curvature}, respectively.

\bigskip
 If $F$ is a Berwald metric, then ${\bf J}=0$ and ${\bf S}=0$.  The converse is true too in dimension two, but it is not clear in higher dimensions (Cf. \cite{Sh2}).

\section{Finsler metrics with constant S-curvature}

The following lemma is cruial for the proof of Theorem \ref{thmnonpositive}. 
 
\begin{lem} Let $(M, F)$ be an $n$-dimensional Finsler manifold. 
Suppose that there is a constant $c$ and a closed $1$-form $\gamma$ such that 
\[
{\bf S}(x,y)= (n+1) c F(x,y)+\gamma_x(y), \ \ \ \ \ y\in T_xM,\]
 then along any geodesic $\sigma=\sigma(t)$, the vector field ${\bf I}(t):= I^i(\sigma(t), \dot{\sigma}(t)) \pxi|_{\sigma(t)}$ satisfies the following equation:
\be
{\rm D}_{\dot{\sigma}}{\rm D}_{\dot{\sigma}}  {\bf I} (t) + {\bf R}_{\dot{\sigma}(t)} ( {\bf I}(t) ) =0.\label{SSKK1}
\ee
\end{lem}
{\it Proof}: It is known that 
the Landsberg tensor  satisfies the following equation \cite{Mo1} \cite{MoSh} :
\be
J_{k|m}y^m + I_mR^m_{\ k}  = -  {1\over 3} \Big \{ 2 R^m_{\ \; k\cdot m} + R^m_{\ \; m\cdot k} \Big \} \label{Moeq2}
\ee
and the 
 S-curvature satisfies
the following equation  \cite{CMS} \cite{Mo2}:
\be
{\bf S}_{\cdot k|m}y^m - {\bf S}_{| k} 
= - {1\over 3} \Big \{ 2 R^m_{\ \; k\cdot m} + R^m_{\ \; m\cdot k} \Big \}.\label{SSKKMo}
\ee
It follows from (\ref{Moeq2}) and (\ref{SSKKMo}) that 
\[
J_{k|m}y^m +I_m R^m_{\ k} = {\bf S}_{\cdot k|m}y^m - {\bf S}_{|k}. \label{EIILJ}
\]
By (\ref{Pg3***}), we can rewrite the above equation as follows
\be
I^i_{\ | p|q} y^py^q +  R^i_{\ m}I^m = g^{ik} \Big \{ {\bf S}_{\cdot k|m}y^m - {\bf S}_{|k} \Big \}. \label{EIILJ*}
\ee
Note that  $F= \sqrt{g_{ij}y^iy^j}$ satisfies
\[  F_{|m} ={g_{ij|m}y^iy^j\over 2F}= 0, \  \ \ \ F_{\cdot k |m}={ g_{ik|m}y^i\over F}=0. \]
Since   $\gamma=\gamma_i dx^i$ is closed, it satisfies
\[ \gamma_{\cdot k|m}y^m - \gamma_{|k}
= \Big \{  {\pa \gamma_k \over \pa x^m} - {\pa \gamma_m \over \pa x^k}     \Big \}y^m =0.\]
We have 
\[ {\bf S}_{\cdot k|m}y^m -{\bf S}_{|k}
= (n+1) c \Big \{ F_{\cdot k|m} y^m - F_{|k}\Big \} + \gamma_{\cdot k |m} y^m -\gamma_{|k}=0.\]
Then  (\ref{EIILJ*}) is reduced to 
\be
I^i_{\ | p|q} y^py^q +  R^i_{\ m}I^m =0. \label{JIK*}
\ee
Since $\sigma$ is a geodesic, we have  
\[ D_{\dot{\sigma}}D_{\dot{\sigma}} {\bf I}(t)=I^i_{\ |p|q} (\sigma(t), \dot{\sigma}(t))\dot{\sigma}^p(t)\dot{\sigma}^q(t) \pxi|_{\sigma(t)}.
\]
Then (\ref{JIK*}) restricted to $\sigma(t)$ gives rise to
 (\ref{SSKK1}).
\qed

\section{Proof of Theorem \ref{thmnonpositive}}

In this section, we are going to prove  a slightly more general version of Theorem \ref{thmnonpositive}.

\begin{thm}\label{thmnonpositive*} 
Let $(M, F)$ be an $n$-dimensional complete Finsler manifold with nonpositive flag curvature  ${\bf K}\leq 0$ and  almost constant S-curvature ${\bf S}= (n+1) c F +\gamma$ ($c=constant$ and $\gamma$ is a closed $1$-form).
Suppose that the mean Cartan torsion grows sub-linearly. Then 
${\bf J}=0$ and ${\bf K}(P, y)=0$ for the flag $P={\rm span}\{ {\bf I}_y, y \}$ whenever ${\bf I}_y \not=0$. Moreover $F$ is Riemannian at points
where ${\bf K}<0$. 
\end{thm}
{\it Proof}: Let $y\in T_xM$ be an arbitrary non-zero vector and let $\sigma=\sigma(t)$ be the geodesic with $\sigma(0)=x$ and $\dot{\sigma}(0)=y$. Since the Finsler metric is complete, one may assume that $\sigma$ is   defined on $(-\infty, \infty)$. 
${\bf I}$ and ${\bf J}$ restricted to $\sigma$ are vector fields along $\sigma$, 
\[ {\bf I}(t):= I^i\Big (\sigma(t), \dot{\sigma}(t)\Big ) \pxi|_{\sigma(t)}, \ \ \ \ \ {\bf J}(t):= J^i\Big (\sigma(t), \dot{\sigma}(t)\Big ) \pxi|_{\sigma(t)}.\]
It follows from (\ref{Pg3***}) that 
\be
D_{\dot{\sigma}}{\bf I}(t)  = I^i_{\ |m}\Big  (\sigma(t), \dot{\sigma}(t)\Big )\dot{\sigma}^m(t) \pxi|_{\sigma(t)}= {\bf J}(t). \label{DIJ}
\ee

If ${\bf I}(t) \equiv 0$, then by (\ref{DIJ}),
${\bf J}_y = D_{\dot{\sigma}} {\bf I}(0)=0$. 
From now on, we assume that  ${\bf I}(t)\not\equiv 0$.
Let
\be \varphi(t):= \sqrt{\g_{\dot{\sigma}(t)} \Big ( {\bf I}(t), {\bf I}(t) \Big )}.\label{varphi}
\ee
Let $I= (a,b)\not=\emptyset$ be a maximal interval on which $\varphi(t) >0$. 
We have
\[
  \varphi \varphi'=  {\bf g}_{\dot{\sigma}} \Big ({\bf I}, D_{\dot{\sigma}}{\bf I } \Big )  \leq \sqrt{{\bf g}_{\dot{\sigma}} \Big ({\bf I}, {\bf I } \Big )}
\sqrt{{\bf g}_{\dot{\sigma}} \Big (D_{\dot{\sigma}}{\bf I}, D_{\dot{\sigma}}{\bf I }  \Big ) }
 =  \varphi \sqrt{{\bf g}_{\dot{\sigma}} \Big (D_{\dot{\sigma}}{\bf I}, D_{\dot{\sigma}}{\bf I } \Big ) }.
\]
This is,
\be
\varphi' \leq  \sqrt{{\bf g}_{\dot{\sigma}} \Big (D_{\dot{\sigma}}{\bf I}, D_{\dot{\sigma}}{\bf I } \Big )}.\label{s11}
\ee
By assumption ${\bf K}\leq 0$ and (\ref{s11}), we have 
\begin{eqnarray}
\frac{1}{2} [ \varphi^2]''& = & \g_{\dot{\sigma}} \Big ({\rm D}_{\dot{\sigma}}{\rm D}_{\dot{\sigma}}  {\bf I} , {\bf I} \Big ) 
+ \g_{\dot{\sigma}} \Big ({\rm D}_{\dot{\sigma}}  {\bf I}, {\rm D}_{\dot{\sigma}}  {\bf I} \Big )\nonumber\\
& = & -\g_{\dot{\sigma}} \Big ({\bf R}_{\dot{\sigma}} ({\bf I}),  {\bf I} \Big ) 
+  \g_{\dot{\sigma}} \Big ({\rm D}_{\dot{\sigma}}  {\bf I} , {\rm D}_{\dot{\sigma}}  {\bf I} \Big )\nonumber\\
& \geq &   \g_{\dot{\sigma}} \Big ({\rm D}_{\dot{\sigma}}  {\bf I} , {\rm D}_{\dot{\sigma}}  {\bf I} \Big )   \geq    \varphi'^2.\label{eqv1}
\end{eqnarray}
We obtain that $\varphi''(t) \geq 0$. 

We claim that $\varphi'(t) \equiv 0$.
Suppose that 
$\varphi'(t_o)\not=0$ for some $t_o \in I$. 
If $\varphi'(t_o) >0$, then  
\[ \varphi(t) \geq \varphi'(t_o) (t-t_o) +\varphi(t_o), \ \ \ \ \ t > t_o.\]
Thus $b=+\infty$.
If $ \varphi'(t_o) <0$, then  
\[ \varphi(t) \geq \varphi'(t_o) (t-t_o) +\varphi(t_o) >\varphi(t_o)>0, \ \ \ \ \ t < t_o.\]
Thus $a=-\infty$.
In either case, $\varphi(t)$ grows at least linearly. Note that for $p=\sigma(t_o)$,
$  {\cal I}_p(| t-t_o|)\geq \varphi(t)$. 
We see that ${\bf I}$ grows at least linearly. 
This is impossible. 
Thus 
$\varphi'(t)\equiv 0$ and $\varphi(t) = constant >0$. In this case, $I= (-\infty, \infty)$.

It follows from (\ref{eqv1}) that 
\[ \g_{\dot{\sigma}} \Big ({\bf R}_{\dot{\sigma}} ({\bf I}),  {\bf I} \Big )=0, \ \ \ \ 
\
{\rm D}_{\dot{\sigma}}  {\bf I}=0.\]
By (\ref{DIJ}), we get ${\bf J}_y=D_{\dot{\sigma}}{\bf I}(0)=0$. Since ${\bf I}_y$ is orthogonal to $y$ with respect to ${\bf g}_y$, 
${\bf K}(P, y)=0$ for 
$P= {\rm span} \{{\bf I}_y, y\}$ whenever ${\bf I}_y \not=0$.

\bigskip

Assume that ${\bf K} <0$ at a point $x\in M$. It follows from ${\bf g}_y ({\bf R}_y ({\bf I}_y) , {\bf I}_y )=0$ that ${\bf I}_y =0$ for all $y\in T_xM\setminus\{0\}$. 
By  Deicke's theorem \cite{De}, $F$ is Riemannian. \qed

\bigskip
Two natural problems arise: 
\ben
\item[(a)] Is there any {\it complete} non-Landsberg metric on $\R^n$ with ${\bf K} \leq 0$,  ${\bf S}= (n+1) cF$ and 
${\cal I}_p (r) \sim C r$ (as $r \to +\infty$)?
\item[(b)] What is the metric structure of a complete Finsler metric on $\R^n$ ($n\geq 3$) satisfying 
${\bf K}= 0$, ${\bf J}=0$ and ${\bf S}=0$ ?
\een

\bigskip
If the flag curvature is strictly negative, we have the following 

\begin{thm}\label{thm4.2}
Let $(M, F)$ be an $n$-dimensional complete Finsler manifold with ${\bf K} \leq -1$ and almost constant S-curvature. Suppose that the mean Cartan torsion ${\bf I}$ grows sub-exponentially a rate of $k=1$. Then $F$ is Riemannian. 
\end{thm}
{\it Proof}: The proof is similar. Assume that ${\bf I}_y \not=0$ for some non-zero vector $y\in T_xM$.
Let $y\in T_xM$ be an arbitrary vector and $\sigma=\sigma(t)$ be the geodesic with $\sigma(0)=x$ and $\dot{\sigma}(0)=y$. 
Let $\varphi(t)$ be defined by (\ref{varphi}). 
Let $I= (a,b)\not=\emptyset$ be the maximal interval on which $\varphi(t) >0$ and $0\in I$.  
By assumption ${\bf K}\leq -1$ and (\ref{s11}), we obtain
\begin{eqnarray*}
\frac{1}{2} [ \varphi^2]''& = & \g_{\dot{\sigma}} \Big ({\rm D}_{\dot{\sigma}}{\rm D}_{\dot{\sigma}}  {\bf I} , {\bf I} \Big ) 
+ \g_{\dot{\sigma}} \Big ({\rm D}_{\dot{\sigma}}  {\bf I}, {\rm D}_{\dot{\sigma}}  {\bf I} \Big )\nonumber\\
& = & -\g_{\dot{\sigma}} \Big ({\bf R}_{\dot{\sigma}} ({\bf I}),  {\bf I} \Big ) 
+  \g_{\dot{\sigma}} \Big ({\rm D}_{\dot{\sigma}}  {\bf I} , {\rm D}_{\dot{\sigma}}  {\bf I} \Big )\\
& \geq & \varphi^2 + \varphi'^2.
\end{eqnarray*}
This gives rise to the following inequality
\be \varphi'' -\varphi \geq 0.\label{eqv}
\ee

We claim that $\varphi'(t) \equiv 0$. 
Suppose that $\varphi'(t_o) \not=0$ for some $t_o\in I$. 
Let 
\[ \varphi_o(t):= \varphi(t_o) \cosh (t-t_o) + \varphi'(t_o) \sinh (t-t_o).\]
Let $h(t):= \varphi'(t)/\varphi(t)$ and  $h_o(t):= \varphi'_o(t)/\varphi_o(t)$.  
\[ \chi(t):=e^{\int_{t_o}^t (h(\tau)+h_o(\tau) ) d\tau } \Big [ h(t)- h_o(t) \Big ].\] 
It is easy to verify that $\chi'(t) \geq 0 $ and $\chi(t_o)=0$. Thus $\chi (t) \geq 0$ for $t > t_o$ and $\chi(t) \leq 0$ for $t < t_o$. This implies that
$ h(t) \geq h_o(t)$ for $t > t_o$ and $ h(t)\leq h_o(t)$ for $t < t_o$. Note that 
\[ h(t)-h_o(t)= \frac{ \varphi'(t)}{\varphi(t)}
- \frac{\varphi'_o(t)}{\varphi_o(t)}
= \frac{d}{dt} \Big [ \ln \frac{\varphi(t)}{\varphi_o(t)} \Big ] .\]
Thus $ [\varphi/\varphi_o]'(t)\geq 0$  for $t >t_o$ and 
$[\varphi/\varphi_o]'(t)\leq  0$ for $t < t_o$. 
We conclude that 
\[ \varphi(t) \geq \varphi_o(t), \ \ \ \ \ \ a < t < b .\]
 If $ \varphi'(t_o) > 0$, then 
\[ \varphi(t) \geq \varphi_o(t) >0, \ \ \ \ \ t > t_o.\]
Thus $b =+\infty$ and 
\[\liminf_{t \to +\infty} \frac{\varphi(t)}{e^{t-t_o}} \geq 
\frac{\varphi(t_o) +\varphi'(t_o) }{2} >0.\] 
If $\varphi'(t_o) < 0$, then 
\[ \varphi(t) \geq \varphi_o(t) >0, \ \ \ \ \ t < t_o.\]
Thus $ a =-\infty$ and 
\[ \liminf_{t\to -\infty} \frac{\varphi(t)}{e^{t-t_o}} \geq 
\frac{\varphi(t_o) - \varphi'(t_o) }{2} >0.\] 
Note that ${\cal I}_p(|t-t_o|) \geq \varphi(t)$ for $p=\sigma(t_o)$. Thus ${\bf I}$ grows at least exponentially at rate of $k=1$.
 This contradicts the assumption. Therefore, $\varphi' (t) \equiv 0$.

Since $\varphi'(t)\equiv 0$, we conclude that  $\varphi(t) \equiv 0$ by (\ref{eqv}). In particular, ${\bf I}_y =\varphi(0) =0$. This contradicts the assumption at the beginning of the argument. 

Therefore  ${\bf I} \equiv 0$ and $F$ is Riemannian by Deicke's theorem \cite{De}. 
\qed

\bigskip
A natural problem arises: Is there any non-Riemannian {\it complete} Finsler metric on 
$\R^n$ satisfying ${\bf K} \leq -1$,  ${\bf S}= (n+1)cF$, but ${\cal I}_p(r)\sim Ce^r $ (as $r \to +\infty$)? This problems remains open.

\bigskip
\begin{ex}\label{exSh}{\rm Let $\phi=\phi(y)$ be a Minkowski norm on $\R^n$ and ${\cal U}
:= \{ y\in \R^n \; | \; \phi(y) <1 \}$. Let $\Theta= \Theta (x, y)$ be a function on $T{\cal U}\cong {\cal U}\times \R^n$ defined by 
\[ \Theta (x, y) = \phi \Big ( y - \Theta(x,y) x \Big ).\]
$\Theta$ is a Finsler metric on ${\cal U}$ which is called the {\it Funk metric} \cite{Fk1}. The Funk metric satisfies the following important equation
\be
\Theta_{x^k} (x,y) =\Theta(x,y)\Theta_{y^k}(x,y).\label{Funkeq}
\ee
Let $a\in \R^n$ be an arbitrary constant vector $a\in \R^n$ with $|a| <1$. Let 
\[ F:= \Theta(x, y) + {\langle a, y\rangle \over 1+\langle a, x \rangle }, \ \ \ \ y\in T{\cal U}\cong {\cal U}\times \R^n.\]
Clearly, $F$ is a Finsler metric near the origin. 
By (\ref{Funkeq}), one sees that 
the spray coefficients of $F$ are given by $G^i= P y^i$, where
\[
P :={1\over 2}  \Big \{ \Theta(x, y) - {\langle a, y \rangle \over 1 + \langle a, x \rangle } \Big \} .
\]
Then using the  above formula for $G^i$ and (\ref{Riemann}), one can easily show that  $F$ has constant flag curvature ${\bf K}= -{1\over 4}$ (see Example 5.3 in \cite{Sh4}). Now let us compute the S-curvature of $F$. 
A direct computation gives 
\[ {\pa G^m\over \pa y^m} = (n+1) P .\]
Let $dV= \sigma_F(x) dx^1 \cdots dx^n$ be the Finsler volume form on $M$. 
Using (\ref{Scurvature}), we obtain 
\begin{eqnarray*}
{\bf S} & = & (n+1) P(x,y) - y^m {\pa \over \pa x^m} \Big ( \ln \sigma_F(x)  \Big )\\
& = & {n+1\over 2} F(x,y) - (n+1) { \langle a, y \rangle \over 1+ \langle a, x \rangle } 
- y^m {\pa \over \pa x^m} \Big ( \ln \sigma_F (x) \Big )\\
& = & {1\over 2} (n+1)  F(x,y) + d\varphi_x(y),
\end{eqnarray*}
where 
\be
 \varphi  := - \ln \Big [ (1+ \langle a, x \rangle )^{n+1}  \sigma_F(x)^{1\over n+1} \Big ].\label{hsigma_F}
\ee
Thus $F$ has almost constant S-curvature. Note that $F$ is not Riemannian in general. 
}
\end{ex}

\bigskip
Example \ref{exSh} shows that the completeness in Theorem \ref{thmnonpositive*} can not be replaced by the positive completeness.

\section{An Example}\label{sectionExample}

The local/global structures of Berwald metrics have been completely determined by Z.I. Szabo \cite{Sz},
 but their curvature properties have not been discussed throughly.  Here we are going to compute the Riemann curvature and the mean Cartan torsion for a special class of Berwald 
manifolds constructed from a pair of Riemannian manifolds. 
Then we show that  these metrics satisfy the conditions and the conclusions in Theorem \ref{thmnonpositive} (b).

Let $(M_i, \alpha_i)$, $i=1, 2$,  be arbitrary Riemannian manifolds and $M=M_1\times M_2$.
Let $f: [0, \infty)\times [0, \infty)\to [0, \infty)$ be an arbitrary  $C^{\infty}$ function satisfying 
\[
f(\lambda s, \lambda t ) = \lambda f(s, t), \ \ (\lambda >0)\ \ \ {\rm and}  \ \ \  f(s, t)\not=0 \ {\rm if} \ (s, t)\not=0.\]
Define
\be
 F:=\sqrt{ f \Big ( [\alpha_1 (x_1, y_1)]^2, \ [\alpha_2 (x_2, y_2)]^2\Big )}, \label{Ff}
\ee
where $ x= (x_1, x_2)\in M$ and $y = y_1\oplus y_2 \in T_{(x_1, x_2)} (M_1\times M_2)\cong T_{x_1}M_1 \oplus T_{x_2}M_2$. 
Clearly, $F$ has the following properties:
\ben
\item[(a)] $F (x, y) \geq 0$ with equality holds if and only if $y=0$;
\item[(b)] $F(x, \lambda y) =\lambda F(x, y)$, $\lambda >0$;
\item[(c)] $F(x,y)$ is $C^{\infty}$ on $TM \setminus\{0\}$. 
\een

Now we are going to find additional condition on $f=f(s,t)$ under which the matrix  $g_{ij}:= {1\over 2} [F^2]_{y^iy^j}$ is positive definite.  
Take standard local coordinate systems $(x^a, y^a)$ in $TM_1$ and 
$(x^{\alpha}, y^{\alpha} )$ in $TM_2$. Then $(x^i, y^j ):= (x^a, x^{\alpha}, y^a, y^{\alpha})$ is a standard local coordinate system in $TM$. 
Express 
\[ \alpha_1 (x_1, y_1)= \sqrt{ \bar{g}_{ab}(x_1) y^ay^b }, \ \ \ \ \ 
\alpha_2(x_2, y_2)  = \sqrt{ \bar{g}_{\alpha\beta}(x_2)y^{\alpha}y^{\beta} },\]
where $y_1 = y^a {\pa \over \pa x^a}$ and $y_2 = y^{\alpha} {\pa \over \pa x^{\alpha}}$.
We obtain
 \be \Big (  g_{ij}  \Big ) = 
\pmatrix{2 f_{ss} \bar{y}_a \bar{y}_b + f_s \bar{g}_{ab} & 2 f_{st} \bar{y}_a \bar{y}_{\beta}  \cr
	2 f_{st} \bar{y}_b \bar{y}_{\alpha} & 2 f_{tt} \bar{y}_{\alpha} \bar{y}_{\beta}+ f_t \bar{g}_{\alpha\beta} \cr},\label{gijgab}
\ee
where $ \bar{y}_a := \bar{g}_{ab} y^b$ and $\bar{y}_{\alpha}:= \bar{g}_{\alpha\beta} y^{\beta}$. By an elementary argument, one can show that $\Big ( g_{ij} \Big )$ is positive definite if and only if  $f$ satisfies the following conditions:
\[ f_s >0, \ \ \ \ f_t > 0, \ \ \ \ f_s + 2 s f_{ss} >0, \ \ \ \ f_t + 2t f_{tt} >0,\]
and 
\[ f_s f_t -2ff_{st} 
 >0.     \]
 In this case, 
\be
\det \Big ( g_{ij} \Big ) 
= h\Big ( [\alpha_1]^2, \ [\alpha_2]^2 \Big ) \det \Big ( \bar{g}_{ab} \Big ) 
\det \Big ( \bar{g}_{\alpha\beta} \Big ),  \label{detBer1}
\ee
where 
\[
h := ( f_s)^{n_1-1} (f_t)^{n_2-1} 
\Big \{ f_s f_t - 2 f f_{st} \Big \}, 
\]
where $n_1:=\dim M_1$ and $n_2:=\dim M_2$.

By a direct computation, one knows that the
 spray coefficients of $F$ are splitted as the direct sum of the spray coefficients of $\alpha_1$ and $\alpha_2$, that is,
\be
  G^a(x,y) = \bar{G}^a(x_1,y_1), \ \ \ \ \ G^{\alpha}(x,y) =\bar{G}^{\alpha}(x_1, y_1),\label{GGGab}
\ee
where $\bar{G}^a$ and $\bar{G}^{\alpha}$ are the spray coefficients of $\alpha_1$ and $\alpha_2$ respectively. 
From (\ref{GGGab}), one can see that the spray of $F$ is independent of the choice of a particular function $f$. In particular, $G^i$ are quadratic in $y\in T_xM$.  Thus $F$ is a Berwald metric. This  fact is claimed in \cite{Sz}. Since $F$ is a Berwald metric,    ${\bf J}=0$ and ${\bf S}=0$ \cite{Sh2}.

The Riemann tensor of $F$ is given by
\[  \Big ( R^i_{\ j} \Big ) = \Big ( \bar{R}^i_{\ j} \Big )= \pmatrix{ 
 \bar{R}^a_{\ b}&  0 \cr 0  & \bar{R}^{\alpha}_{\ \beta} \cr} ,\]
where $\bar{R}^a_{\ b}$ and $\bar{R}^{\alpha}_{\ \beta} $ are the coefficients of the Riemann tensor of $\alpha_1$ and $\alpha_2$ respectively. 
Let $R_{ij}:= g_{ik}R^k_{\ j}$, $\bar{R}_{ab}:= \bar{g}_{ac} \bar{R}^c_{\ b}$ and $\bar{R}_{\alpha\beta} := \bar{g}_{\alpha\gamma} \bar{R}^{\gamma}_{\ \beta}$. 
Using (\ref{gijgab}), one obtains 
\[ \Big ( R_{ij} \Big ) = \pmatrix{ f_s \bar{R}_{ab} & 0 \cr 
0 & f_t \bar{R}_{\alpha\beta} \cr }.\]
For any vector $v =v^i \pxi|_x \in T_xM$, 
\be
 \g_y \Big ( {\bf R}_y(v), v\Big )
= f_s \bar{R}_{ab}v^av^b + f_t \bar{R}_{\alpha\beta} v^{\alpha} v^{\beta} .\label{Ksign}
\ee
It follows from (\ref{Ksign}) that if $\alpha_1$ and $\alpha_2$ both have nonpositive sectional curvature, then $F$ has nonpositive flag curvature. 

Using (\ref{detBer1}), one can compute the mean Cartan torsion. First, observe that
\[ I_i =  {\pa \over \pa y^i} \Big [ \ln \sqrt{ \det (g_{jk}) } \Big ]\nonumber\\
= {\pa \over \pa y^i} \Big [ \ln \sqrt{ h \Big ( [\alpha_1]^2, [\alpha_2]^2 \Big ) } \Big ].
\]
One obtains
\[
I_a= { h_s \over h} \bar{y}_a \ \ \ \ \ 
I_{\alpha}={h_t \over h} \bar{y}_{\alpha},\]
where $\bar{y}_a:= \bar{g}_{ab}y^b$ and $\bar{y}_{\alpha} := \bar{g}_{\alpha\beta} y^{\beta}$.
Since $ \bar{y}_a \bar{R}^a_{\ b}=0$ and $ \bar{y}_{\alpha} \bar{R}^{\alpha}_{\ \beta}=0$, one obtains 
\[\g_y \Big ( {\bf R}_y ( {\bf I}_y), {\bf I}_y \Big )  = I_i R^i_{\ j} I^j 
= {h_s\over h} \bar{y}_a \bar{R}^a_{\ b} I^b  + {h_t\over h} \bar{y}_{\alpha} \bar{R}^{\alpha}_{\ \beta}I^{\beta} =0.\]
Therefore $F$ satisfies the conditions and conclusions in Theorem \ref{thmnonpositive*}.

\noindent
Zhongmin Shen\\
Dept of Math, IUPUI\\
 402 N. Blackford Street\\
Indianapolis, IN 46202-3216, USA.  \\
zshen@math.iupui.edu\\
www.math.iupui.edu/\~{}zshen


\begin{thebibliography}{Ma}


\bibitem[1]{AZ} H. Akbar-Zadeh, {\it Sur les espaces de Finsler \'{a} courbures sectionnelles constantes}, Bull. Acad. Roy. Bel. Cl, Sci, 5e S\'{e}rie - Tome LXXXIV (1988), 281-322.

\bibitem[2]{BaChSh2}
D. Bao, S. S. Chern and Z. Shen,
{\it Rigidity issues on Finsler surfaces}, Rev. Roumaine Math. Pures Appl. 
 {\bf 42}(1997), 707-735.


\bibitem[3]{BaRo} D. Bao and C. Robles, {\it On Randers metrics of constant curvature}, Reports on Mathematical Physics,  (to appear).

\bibitem[4]{BaRoSh} D. Bao, C. Robles and Z. Shen, {\it Zermelo Navigation on Riemannian manifolds}, preprint. 


\bibitem[5]{BaSh} D. Bao and Z. Shen, {\it Finsler metrics of constant curvature on the Lie group $S^3$}, J. of the  London Math. Soc. {\bf 66}(2002), 453-467.





\bibitem[6]{Be1} 
L. Berwald, {\it Untersuchung der Kr\"{u}mmung allgemeiner metrischer R\"{a}ume auf Grund des in ihnen herrschenden Parallelismus}, Math. Z.
{\bf 25}(1926), 40-73.

\bibitem[7]{Be2} 
L. Berwald, {\it Parallel\"{u}bertragung in allgemeinen R\"{a}umen},
Atti Congr. Intern. Mat. Bologna  {\bf 4}(1928), 263-270.


\bibitem[8]{CMS} X. Chen, X. Mo and Z. Shen, {\it On the flag curvature of  Finsler metrics
of scalar curvature},  J. of the London Math. Soc., (to appear). 


\bibitem[9]{ChSh}
X. Chen and Z. Shen, 
{\it Randers metrics with special curvature properties}, 
Osaka J. of Math. {\bf 40}(2003), 87-101.



\bibitem[10]{Ch}
 S. S. Chern, 
{\it  On the Euclidean connections in a Finsler space},  
Proc. National Acad. Soc., {\bf 29}(1943),  
33-37; or Selected Papers, vol. I\negthinspace I, 107-111, Springer 
1989.

\bibitem[11]{De} A. Deicke, {\it \"{U}ber die Finsler-Raume mit $A_i=0$}, Arch. Math. {\bf 4}(1953), 45-51.


\bibitem[12]{Fk1} P. Funk, {\it \"{U}ber Geometrien, bei denen die Geraden die K\"{u}rzesten sind}, Math. Annalen {\bf 101}(1929), 226-237.



\bibitem[13]{Mo1} X. Mo, {\it The flag curvature tensor on a closed Finsler space}, Results in Math. {\bf 36}(1999), 149-159. 


\bibitem[14]{Mo2} X. Mo,
{\it On the flag curvature of a Finsler space with constant S-curvature}, Houston J. of Math. (to appear).

\bibitem[15]{MoSh} X. Mo and Z. Shen,
{\it On negatively curved Finsler manifolds of scalar curvature}, Canadian Mathematical Bulletin, to appear.

\bibitem[16]{Sh1} Z. Shen,
{\it Volume comparison  and its applications in Riemann-Finsler
geometry}, Advances in Math.
{\bf 128}(1997), 306-328.



\bibitem[17]{Sh2} Z. Shen, {\it Lectures on Finsler Geometry}, World Scientific, Singapore (2001), 307 pages. 

\bibitem[18]{Sh3} Z. Shen, {\it Projectively flat Randers metrics with constant flag curvature}, Math. Ann. {\bf 325}(2003), 19-30.

\bibitem[19]{Sh4} Z. Shen, {\it Projectively flat Finsler metrics of constant flag curvature},  Trans. of Amer. Math. Soc. 
{\bf 355}(4)(2003), 1713-1728. 
\bibitem[20]{Sh5} Z. Shen, {\it Finsler metrics with ${\bf K}=0$ and ${\bf S}=0$}, Canadian J. Math. {\bf 55}(2003), 112-132.

\bibitem[21]{Sz} 
Z.I. Szab\'{o}, 
{\it Positive definite Berwald spaces (Structure theorems on Berwald spaces)}, 
 Tensor,  N. S.
{\bf  35}(1981),  25--39. 



\end{thebibliography}
\end{document}